\documentstyle[12pt]{article}

\def\cl{\centerline}

\begin{document}

\centerline{\bf QUASI-EXACTLY SOLUBLE POTENTIALS AND DEFORMED OSCILLATORS}

\bigskip\bigskip\bigskip
\cl{Dennis BONATSOS}

\cl{European Centre for Theoretical Studies in Nuclear 
Physics and Related Areas (ECT$^*$)}

\cl{Strada delle Tabarelle 286, I-38050 Villazzano (Trento), Italy}

\cl{and}

\cl{Institute of Nuclear Physics, NCSR ``Demokritos''}

\cl{GR-15310 Aghia Paraskevi, Attiki, Greece}

\cl{C. DASKALOYANNIS}

\cl{Department of Physics, Aristotle University of Thessaloniki}

\cl{GR-54006 Thessaloniki, Greece}

\cl{Harry A. MAVROMATIS}

\cl{Physics Department, King Fahd University of Petroleum and 
Minerals}

\cl{ Dhahran 31261, Saudi Arabia}

\bigskip\bigskip\bigskip
\cl{\bf Abstract}

It is proved that quasi-exactly soluble potentials corresponding to 
an oscillator with harmonic, quartic and sextic terms, 
for which the $n+1$ lowest levels of a given parity can be  determined 
exactly, may be approximated by  WKB equivalent potentials corresponding to 
deformed anharmonic oscillators of SU$_q$(1,1) symmetry, which have been used 
for the description of vibrational spectra of diatomic molecules. 
This connection allows for the immediate approximate determination of the 
levels of the same parity lying above the lowest $n+1$ known levels, 
as well as of all levels of the opposite parity. Such connections are 
not possible in the cases of the q-deformed oscillator, the Q-deformed 
oscillator, and the modified P\"oschl-Teller potential with SU(1,1) symmetry.

\bigskip\bigskip\bigskip
\noindent

Talk given by D. Bonatsos at the {\it 6th Hellenic Symposium on Nuclear 
Physics} (26--27 May 1995, Piraeus, Greece). To be published in the 
Proceedings, edited by C. N. Panos.

\vfill\eject

\section{Introduction}

Quantum algebras (also called quantum groups), which are deformed
 versions of the usual Lie algebras, to which they reduce when the deformation 
parameter is set equal to unity,  have  recently been attracting considerable 
attention. The interest in possible physical applications was triggered by the 
introduction of the $q$-deformed harmonic oscillator by Biedenharn and 
Macfarlane in 1989, although similar oscillators had already been in 
existence. The common mathematical structure of these various types 
of oscillators was shown by several authors (see \cite{PLB307} for a list 
of references), using the concept of the generalized deformed oscillator. 

One way to clarify the physical content of deformed oscillators is to try
to construct potentials which give, exactly or  approximately, the same 
spectrum as these oscillators. Along these lines, WKB equivalent potentials 
corresponding to the $q$-deformed oscillator, as well as to deformed 
anharmonic oscillators found useful in the description of vibrational 
spectra of diatomic molecules \cite{BAR} have been constructed, both 
numerically and in analytic form \cite{JMP33}. These 
potentials give approximately the same spectrum as the deformed oscillators 
under discussion. 

On the other hand, several quasi-exactly soluble potentials (QESPs) have been 
introduced \cite{TurUsh,Tur118,Tur276}, 
for which the exact calculation of the first $n+1$ 
energy levels is possible, with no information provided for the rest of the 
levels. In particular, in the case of the harmonic oscillator with both 
quartic and sextic 
anharmonicities \cite{TurUsh}, the exact calculation of the first 
$n+1$ energy levels {\sl of a given parity}  is possible.   

It is therefore interesting to check if there is any relation between the 
WKB equivalent potentials (WKB-EPs) approximating the behavior of the deformed 
oscillators on the one hand and the quasi-exactly soluble potentials on the
other. If such a link exists, as we shall demonstrate in this letter, 
one can fix the parameters of a WKB potential 
equivalent to an appropriate oscillator so that the potential agrees, 
up to the sextic term, with a given quasi-exactly soluble potential. 
If the approximations involved (WKB method, truncation of the WKB 
potential to sixth order only) are good, one can subsequently  use
the levels of the relevant oscillator in order to approximate the 
unknown higher levels of the quasi-exactly soluble potential. In the case
of the harmonic oscillator with both quartic and sextic anharmonicities 
\cite{TurUsh}, the levels of the opposite parity can also be approximated in 
this way.   

\section{Deformed oscillators}

The q-numbers are defined as
$$[x]= {q^x -q^{-x}\over q-q^{-1}}.\eqno(1)$$
In the case that $q=e^{\tau}$, where $\tau$ real, they can be
written as
$$[x]={\sinh (\tau x) \over \sinh (\tau)}, \eqno(2)$$
while in the case that $q$ is a phase ($q=e^{i\tau}$, with $\tau$
real) they take the form
$$[x]={\sin (\tau x)\over \sin (\tau)}. \eqno(3)$$
It is clear that in both cases the q-numbers reduce to the usual
numbers as $q\rightarrow 1$ ($\tau\rightarrow 0$).
 
The q-deformed harmonic oscillator is determined by 
 the creation and annihilation operators $a^+$ and $a$, which
satisfy the  relations:
$$ a a^+ -q^{\mp 1} a^+ a = q^{\pm N}, \quad
[N, a^+] = a^+, \quad [N, a]= -a , \eqno(4)$$
where $N$ is the number operator.
The Hamiltonian of the q-deformed oscillator is
$$ H= {\hbar \omega \over 2} (a a^+ + a^+ a), \eqno(5)$$
and its   eigenvalues are
$$ E(n) = {\hbar \omega\over 2} ([n]+[n+1]). \eqno(6)$$

From Eq. (6) it is clear that the spectrum of the q-deformed
harmonic oscillator is not equidistant. For real $q$ the
spectrum increases more rapidly than in the classical case,
while for $q$ complex it increases less rapidly than
in the classical case, i.e. it gets squeezed. 

In addition to the q-numbers, Q-numbers have been defined in the 
literature as 
$$ [x]_Q = {Q^x-1 \over Q-1}.\eqno(7)$$
The corresponding Q-deformed harmonic oscillator is determined by the
creation and annihilation operators $b^+$ and $b$, satisfying the  
commutation relations 
$$  b b^+ - Q b^+ b =1, \quad [N, b^+] = b^+, \quad [N,b]=-b. \eqno(8)$$
The Hamiltonian of the relevant oscillator is 
$$ H = {\hbar \omega \over 2} (b b^+ + b^+ b),\eqno(9)$$
with corresponding eigenvalues
$$ E(n)= {\hbar \omega \over 2} ([n]_Q +[n+1]_Q).\eqno(10)$$
If $Q=e^T$, with $T$ real, the spectrum increases more rapidly than an
equidistant spectrum for $T>0$, while it increases less rapidly than an 
equidistant spectrum for $T<0$. 

Using two q-deformed oscillators described by the operators
$a_1,a^+_1,N_1$, and $a_2,a^+_2,N_2$, the generators of the
quantum algebra SU$_q$(1,1) are written as
$$K_+=a_1^+a_2^+, \quad K_=a_1a_2, \quad
K_0={1\over 2}\left( N_1+N_2+1\right), \eqno(11)$$
and satisfy the commutation relations
$$ \left[ K_0,K_\pm \right] = \pm K_\pm, \quad
\left[ K_+,K_- \right] = - \left[ 2K_0 \right],\eqno(12)$$
where the square bracket in the rhs of the last equation is a q-number 
as defined in eq. (1). In ref. \cite{BAR}
this symmetry has been used for the description of
vibrational molecular spectra. The q-deformed anharmonic oscillator
used there is a q-generalization of the anharmonic oscillators used
in the usual Lie algebraic approach to molecular spectroscopy.
(See ref. \cite{BAR} for detailed references.) 

\section{WKB equivalent potentials}

In the case in which $q=e^{i\tau}$, the spectrum of the q-oscillator
(Eq. (6)) can be written as
$$E_n = {\hbar\omega\over 2} {\sin(\tau(n+1/2))\over
\sin (\tau/2)}. \eqno(13)$$
In this case one can see \cite{JMP33}
 that the WKB equivalent potential is given by
$$ V(x) =  \left({ \tau \over 2 \sin (\tau /2) }\right)^2
{ m\omega^2\over 2} x^2$$
$$ \Big[ 1-{8\over 15} \left({x \over { 2 R_e}}\right)^4
+  {4448\over 1575} \left({x \over { 2 R_e}}\right)^8
- {345344\over 675675}\left({x \over { 2 R_e}}\right)^{12}
+\ldots\Big] ,\eqno(14)$$
with 
$$ R_e = {1\over \tau} \left( {\hbar^2\over 2m}\right)^{1/2} \left( 
2\sin (\tau/2) \over \hbar \omega\right)^{1/2} .\eqno(15)$$
For $\tau =0$ the usual harmonic oscillator potential is
obtained.  For $x\to \infty$ the potential
goes to a finite limiting value. 

For $q=e^{\tau}$ with $\tau$ real, the
 q-oscillator has a spectrum given by
$$ E_n = {\hbar \omega\over 2} {{\sinh(\tau(n+1/2))} \over
{\sinh (\tau/2)}}, \eqno(16)$$
while the WKB equivalent potential is given by \cite{JMP33}
$$ V(x) = \left({ \tau\over 2 \sinh (\tau /2) }\right)^2
{ m\omega^2\over 2} x^2 $$
$$ \Big[1+ {8\over 15} \left({x \over { 2 R_h}}\right)^4
+  {4448\over 1575} \left({x \over { 2 R_h}}\right)^8
+ {345344\over 675675}\left({x \over { 2 R_h}}\right)^{12}
+\ldots \Big], \eqno(17)$$
with 
$$R_h = {1\over \tau} \left( {\hbar^2 \over 2 m} \right)^{1/2} \left( 
2 \sinh (\tau/2) \over \hbar \omega \right)^{1/2}.\eqno(18)$$
This WKB-EP gives the classical harmonic oscillator potential
for $\tau=0$, while it  goes to infinity for $x\to \infty$.

Using the same inversion technique one finds that the WKB equivalent potential 
for the Q-deformed harmonic oscillator takes the form \cite{Ioann}
$$ V(x) = V_{min} + {(\ln Q)^2 \over Q} \left( {Q+1\over Q-1}\right) ^2 
{1\over 2} m \omega^2 x^2 $$ 
$$\left[ 1 -{2\over 3} \left( {x\over R'}\right)^2
+ {23\over 45} \left( {x\over R'}\right)^4 -{134\over 315} \left( {x\over R'}
\right)^6 + {5297\over 14172} \left( {x\over R'}\right)^8 -\ldots \right],
\eqno(19)$$
where 
$$V_{min} = {\hbar \omega (\sqrt{Q}-1) \over 2 \sqrt{Q} (\sqrt{Q}+1)},\eqno(20)
$$
and
$$ R'= \left( {\hbar \sqrt{Q} (Q-1) \over \omega m (Q+1) }\right)^{1/2} 
{(\ln Q)^{-1} \over \sqrt{2}}.\eqno(21) $$

We now turn our attention to the anharmonic oscillator with
SU$_q$(1,1) symmetry, which has been found useful in the description
of vibrational molecular spectra \cite{BAR}.
The energy  spectrum of this oscillator, in the case of complex $q$, is
given by \cite{BAR}
$$ E_n = E_0' -A {{\sin(\tau(n-N/2))\sin(\tau(n+1-N/2))} \over
{\sin^2 (\tau)}} , \eqno(22)$$
where
$$N=2n_{max} \quad\quad \hbox{or} \quad\quad N=2n_{max}+1,\eqno(23)$$
with $n_{max}$ corresponding to the last level before the
dissociation limit.
For $\tau \to 0$ a Morse or modified P\"oschl-Teller spectrum 
is obtained.
The WKB equivalent potential in this case is 
$$ V(x)= V_{min} +
{A\over 4} \left({ \tau \sin (N\tau )\over  \sin^2 \tau} \right)^2
u^2$$
$$ \Big[1 - {2\over 3} {\tau^2 \cos (N\tau )\over \sin^2 \tau}u^2
+{1\over 45} { \tau^4 (23\cos^2(N\tau)-6)\over
\sin^4\tau}u^4 $$
$$- {2\over 315}
{ \tau^6 ( 67\cos^2(N\tau)-36)\cos(N\tau) \over \sin^6\tau}u^6+\cdots
\Big], \eqno(24)$$
where
$$ V_{min} = E_0' - A {\cos (\tau) -\cos (\tau N)\over 2 \sin ^2 \tau},\eqno(25)
$$
and 
$$u={ \sqrt{2mA}x \over \hbar}. \eqno(26)$$
It can be easily seen that for $\tau \to 0$ Eq. (24) reduces to the
Taylor expansion of the modified P\"oschl-Teller potential 
$$ V(x) = V_{min} + {A\over 4} u^2 \left[ 1-{2\over 3} u^2 +{17\over 45} u^4
-{62\over 315} u^6 + \ldots \right], \eqno(27)$$
where $$ u= {\sqrt{ 2 m A} \over \hbar} x, \eqno(28)$$
while  the modified P\" oschl--Teller potential in closed form is 
$$ V(x)=V_{min} + {A N^2 \over 4}
\tanh^2\left({ \sqrt{2mA}x \over \hbar}\right).\eqno(29)$$
 This means that the WKB-EP of a system with
$SU_q(1,1)$ symmetry is a deformation of the modified P\"oschl--Teller 
potential.

\section{Relation to Quasi Exactly Soluble Potentials}

On the other hand, it is known \cite{TurUsh,Tur118} that the potentials
$$ V(x) = 8 a^2 x^6 + 8 a b x^4 + 2 [b^2 -(2k+3) a ] x^2, \eqno(30)  $$ 
with $k=2n+r$, where $n=0$, 1, 2, \dots and $r=0$, 1 are quasi-exactly 
soluble. The meaning of this term is the following: for these potentials 
one can construct exactly the first $n+1$ levels with parity $(-1)^r$. 
(The extra factors of 2 in eq. (30) in comparison to refs \cite{TurUsh,Tur118}
are due to the fact that we use the usual form of the Schr\"odinger operator
$$ H= - {\hbar^2 \over 2 m} {d^2 \over d x^2 } + V(x),\eqno(31)$$
while in refs \cite{TurUsh,Tur118} the form 
$$ H = -{\hbar^2 \over m} {d^2 \over d x^2} + V(x) \eqno(32)$$
is used.) 

\subsection{The SU$_q$(1,1) anharmonic oscillator}

We wish to check if the WKB-EPs for the deformed oscillators mentioned 
above correspond, up to terms of $x^6$, to quasi-exactly soluble potentials
of the form of eq. (30). For reasons of convenience we start the comparison
from the SU$_q$(1,1) potential of eq. (24), setting $\hbar = m=1$. 
Comparing the coefficients of $x^2$, $x^4$ and $x^6$ in eqs (24) and (30) 
we get respectively the following equations
$$ 2 ( b^2 -(2k+3) a) = {A^2\over 2} {\tau^2 \sin^2(N\tau) \over 
\sin^4\tau}, \eqno(33)$$ 
$$ 8 a b = -{2\over 3} A^3 {\tau^4 \sin^2(N\tau) \cos(N\tau) \over 
\sin^6\tau}, \eqno(34)$$
$$ 8 a^2 = {2\over 45} A^4 {\tau^6 \sin^2(N\tau) (23\cos^2(N\tau)-6)\over
\sin^8 \tau}.\eqno(35) $$
We can consider $a^2=1$ without loss of generality. 
For the ground state to be normalizable one should have $a\geq 0$ 
\cite{Tur276}, which implies in the present case $a=1$. Having in mind that 
in eq. (22) one must have $A>0$ in order to get an increasing spectrum 
(this fact can be seen also in the realistic applications to molecular 
spectra given in ref. \cite{BAR}), from eq. (35) one finds
$$ A= \sqrt{6\sqrt{5}} {\sin^2\tau \over \tau^{3/2} \sqrt{\sin(N\tau)} 
\root 4 \of {23\cos^2(N\tau)-6} } .\eqno(36) $$
Because of the symmetry $q\leftrightarrow q^{-1}$ characterizing the 
q-deformed oscillators, it suffices to consider $\tau>0$. The following 
conditions should then be satisfied 
$$ \sin(N\tau) >0, \eqno(37)$$
$$ \cos^2(N\tau)>6/23.\eqno(38)$$
Then eq. (34) gives 
$$ b = -\sqrt{ {15\sqrt{5}\over 2} } {\sqrt{\sin(N\tau)} \cos(N\tau) \over
\sqrt{\tau} (23\cos^2(N\tau)-6)^{3/4} },\eqno(39) $$
while eq. (33) gives
$$ 2k+3 = - {3\sqrt{5} \over 2} {\sin(N\tau) (18\cos^2(N\tau) -6) \over 
\tau (23\cos^2(N\tau)-6)^{3/2} }.\eqno(40)$$
The fact that $k=2n+r$ has to be positive gives the additional condition 
$$ \cos^2(N\tau) < 1/3.\eqno(41) $$
From eq. (39) it is clear that $b$ can be either positive (when 
$\cos(N\tau)<0$) or negative (when $\cos(N\tau)>0$). 

For reasons of convenience, let us consider the case $b>0$ first. 
In this case the conditions $\sin(N\tau)>0$ and $\cos(N\tau)<0$ imply that 
we can limit ourselves to the region $\pi > N\tau >\pi/2$. The conditions
(38) and (41) then imply  that $ 2.1863 + 2 \pi l > N\tau > 2.1069 + 2 \pi l$, 
with $l=0$, 1, 2, \dots. Thus though we can find an infinity of values for
$N$, as we shall see it suffices to consider $l=0$.  

We first wish to check the accuracy of our approach, which already involves
two major approximations: the WKB approximation and the omission in the 
WKB-EP of terms higher than $x^6$. In order to achieve this, we will compare
the results given by this method to the exact results given in 
ref. \cite{TurUsh}
 for the case $n=1$, $r=0$, in which $2k+3=7$. From eq. (40) one 
then sees that a possible solution is $N=151$, $\tau=0.0144503$, which 
gives $A=0.4343473$, $b=12.589097$, while from eq. (25) (for $V_{min}=0$) one
 has $E'_0= 1636.8943$. 
We see therefore that eqs (33)--(35) are indeed satisfied and the potentials
of eqs (24) and (30) (up to the $x^6$ term) are identical, being 
$$ V(x) = 303.02 ( x^2 + 0.3324 x^4 + 0.02640 x^6).\eqno(42) $$
The two solutions which can be obtained exactly are \cite{TurUsh} 
$$E_0 = 3 b -2\sqrt{b^2+2} = 12.4307, \qquad E_2 = 3b+ 2\sqrt{b^2+2} = 
63.1039,$$
$$ E_2-E_0 = 50.6731,\eqno(43) $$
while eq. (22) gives the complete spectrum, including 
$$E_0 = 12.3805, \qquad E_2 = 63.0584, \qquad E_2-E_0= 50.6779. \eqno(44) $$
We remark that the agreement between the exact results of ref. 
\cite{TurUsh} and the 
predictions of eq. (22) is excellent, implying that in this case:

i) the WKB method is accurate, 

ii) the omission of terms higher than $x^6$ is a good approximation.

An additional check for the case $n=3$, $r=0$ is also made. In this case 
$2k+3=15$. One possible solution of eqs (33)-(35) in this case is given by 
$N=325$, $\tau=0.00671384$, $A=0.2960795$, $b= 18.469158$, $E'_0=5168.941$. 
The potential is 
$$V(x) = 652.20 (x^2 + 0.2265 x^4 + 0.01227 x^6).\eqno(45) $$ 
Eq. (22) gives the levels
$$ E_0= 18.1126, \qquad E_2 = 91.3482, \qquad E_4=165.8771, \qquad 
E_6 = 241.6456,\eqno(46)$$
while, following the procedure of ref. \cite{TurUsh}, 
one sees that in this case
the exact energy levels are the roots of the equation
$$ E^4 -28 b E^3 + (254 b^2 -240) E^2 + (-812 b^3 + 2592 b) E +
(585 b^4 -4656 b^2 +2880) =0,\eqno(47) $$
given by 
$$ E_0=18.1429, \qquad E_2=91.3783, \qquad E_4=165.913, \qquad E_6=241.703. 
  \eqno(48)$$ 
We remark that excellent agreement between the approximate and the exact 
results is again obtained. 

The agreement remains equally good if more levels are considered. 
For the  case $n=9$, $r=0$, in which $2k+3=39$, one possible solution 
of eqs (33)-(35) is given by $N=399$, $\tau=0.00545864$, $A=0.2703882$, 
$b=21.275801$, $E_0'= 7126.0336$. The potential is 
$$ V(x)= 827.43 (x^2 + 0.2057 x^4 + 0.009669 x^6).$$
The energy levels are compared to the exact ones in Table I.   


\begin{table}
\caption{
Exact energy levels of the quasi-exactly soluble potential of eq. (30)
with $n=9$, $r=0$, $a=1$, $b=21.275801$, compared to the approximate energy 
levels of the WKB-EP potential of eq. (24) with $A=0.2703882$, $E_0'= 
7126.0336$, 
$\tau= 0.00545864$, $N=399$, given by eq. (22). }  
\medskip
\centerline{ 
\begin{tabular}{| r r r |}
\hline\hline
 $n$  & $E_n$ (approx) & $E_n$ (exact)  \\
\hline
0  & 20.3991 & 20.4153 \\
2  & 102.672 & 102.682 \\
4  & 186.131 & 186.138 \\
6  & 270.735 & 270.749 \\
8  & 356.444 & 356.485 \\
10 & 443.217 & 443.317 \\
12 & 531.013 & 531.218 \\
14 & 619.791 & 620.165 \\
16 & 709.507 & 710.133 \\
18 & 800.118 & 801.101 \\
\hline\hline
\end{tabular}  }
\end{table}

We turn now to the case with $b<0$. The conditions $\sin(N\tau)>0$ and 
$\cos(N\tau)>0$ imply that we should limit ourselves to the region 
$\pi/2 > N\tau >0$. Then the conditions (38) and (41) imply that
$ 1.0347 > N\tau > 0.9553$. (We again ignore terms of $2 \pi l$.)  
 
We again consider the case $n=1$, $r=0$, $2k+3=7$. One possible solution is 
given by $N=61$, $\tau=0.0157377$, $A=0.4538508$, $b=-12.108743$, $E'_0= 
390.66689$.
The relevant potential is 
$$ V(x) = 279.095 (x^2 - 0.3471 x^4 + 0.02866 x^6).\eqno(49) $$
The exact solutions given by ref. \cite{TurUsh} are 
$$ E_0 = -60.7083, \qquad E_2 = -11.9441,\eqno(50)$$
while eq. (22) gives 
$$ E_0 = 11.7456, \qquad E_2= 57.3764.\eqno(51) $$ 

The reasons for this dramatic failure are well understood. The potential 
of eq. (49) has a central well in the middle, having its minimum at 
$x=0$, $E=0$, plus two symmetrically located  wells, on either side of the 
central one, with their minima at $E<0$. The method of ref. \cite{TurUsh}
gives the lowest two energy eigenvalues of the system, which in this case are
the lowest two levels in the side wells. The WKB method used in ref. 
\cite{JMP33}, 
however, is known to be valid only for small values of $x$, i.e. in the 
area of the central well. Therefore in this case eq. (51) gives the lowest
two levels of the central well, which in this case are {\sl not} the lowest 
two levels of the system.  

We therefore conclude that a correspondence between the SU$_q$(1,1) WKB-EP
and the quasi-exactly soluble potentials of eq. (30) can be made only 
for $b>0$, in which case both methods will give the levels corresponding 
to a well around $x=0$. In this case, as we have already seen, the 
approximation is very good, providing us with the following method:

 Given a quasi-exactly soluble potential for which only the $n+1$ lowest 
levels of parity $r$ can be found exactly, we can choose appropriately the 
parameters of an SU$_q$(1,1) WKB-EP, so that the two potentials are 
identical up to order $x^6$. Then using eq. (22), i.e. the known eigenvalues 
of  the WKB-EP, we can approximate the levels of the same parity lying above 
the first $n+1$ known levels, as well as the levels of opposite parity, 
for which the method of ref. \cite{TurUsh} gives no information {\sl for the 
potential under discussion}. (It is worth noting that changing the 
values of $n$ and $r$ changes the potential.)  

\subsection{The q-deformed oscillator}

It is now worth examining if the WKB-EPs corresponding to the q-oscillator 
and the Q-oscillator can be made to correspond to a quasi-exactly 
soluble potential.  

In the case of the q-deformed harmonic oscillator with $q=e^{i\tau}$ ($\tau$ 
real) one should set $b=0$, since no $x^4$ terms appear in the WKB-EP of 
eq. (14). Comparison of the coefficients of $x^6$ gives then the condition
$$ 8 a^2 = - {\tau^6 \omega^4 \over 240 \sin^4(\tau/2)} ,\eqno(52)$$
which cannot be satisfied. Therefore no connection between this oscillator and
the quasi-exactly soluble potentials of eq. (30) is possible. 

In the case of the q-deformed harmonic oscillator with $q=e^{\tau}$ ($\tau$ 
real) one should again set $b=0$, since no $x^4$ terms appear in the 
WKB-EP of eq. (17). In addition, because of the symmetry $q\leftrightarrow
q^{-1}$ characterizing this oscillator, it suffices to consider $\tau>0$.
Then, taking into account that $a>0$ for the ground state to be normalizable
\cite{Tur276},  the coefficients of $x^6$ give the equation 
$$ a = {\tau^3 \omega^2 \over 8\sqrt{30} \sinh^2(\tau/2)}, \eqno(53)$$
while   the coefficients of $x^2$ give the condition 
$$ 2k+3 = -\sqrt{15\over 2} {1\over \tau}.\eqno(54)$$
Since both $2k+3$ and $\tau$ are positive, this condition cannot be satisfied. 
Therefore no connection between this oscillator and the quasi-exactly soluble 
potentials of eq. (30) is possible, either. 

\subsection{The Q-deformed oscillator}

In the case of the Q-deformed oscillator the $x^4$ term is present in the 
WKB-EP (see eq. (19)), in contrast to the q-deformed oscillator. In this case 
the coefficients of the $x^6$, $x^4$ and $x^2$ terms give the equations
$$ a=\pm {1\over 2} \sqrt{23\over 45} {(\ln Q)^3 \over Q} \left( Q+1 \over Q-1
\right)^2 \omega^2,\eqno(55)$$
$$b= \mp {1\over 2} \sqrt{5\over 23} {\ln Q \over \sqrt{Q}} {Q+1\over Q-1} 
\omega,\eqno(56)$$
$$ 2k+3 = \mp {27\over 23} \sqrt{5\over 23} {1\over \ln Q}.\eqno(57)$$
We know that for the ground state to be normalizable one should have $a>0$. 
Therefore in eq. (55) the selection of the upper (lower) sign requires 
$\ln Q >0$ ($\ln Q <0$). In both cases it is impossible to satisfy eq. (57). 
Therefore no connection between the Q-deformed oscillator and the quasi-exactly
soluble potentials of eq. (30) is possible. 

\subsection{The modified P\"oschl--Teller potential}

In the case of the modified P\"oschl-Teller potential, comparing its Taylor
expansion (eq. (27)) to the quasi-exactly soluble potentials of eq. (30) and 
taking into account that $a>0$ one obtains
from the coefficients of the $x^6$, $x^4$ and $x^2$ terms the following 
equations
$$ a =\sqrt{17\over 5} {A^2\over 6}, \eqno(58)$$
$$ b =-\sqrt{5\over 17} {A\over 2}, \eqno(59)$$
$$ 2k+3 = -{18\over 17} \sqrt{5\over 17}.\eqno(60)$$  
Since $2k+3$ has to be positive, it is clear from the last equation that no 
connection between the modified P\"oschel-Teller potential and the 
quasi-exactly soluble potentials of eq. (30) is possible. 

\section{Discussion}

The following comments on the results can now be made:

i) The quasi-exactly soluble potentials are known to be related to the SL(2)
\cite{Tur118} and deformed SU(2) \cite{SunLi} symmetries. Among the several 
oscillators 
considered here, the only one for which the connection to QESPs was possible
is the oscillator with SU$_q$(1,1) symmetry. Apparently the H$_q$(4) symmetry
of the q-deformed oscillator is ``not  enough'' to guarantee such a 
connection. 

ii) The QESPs of eq. (30) are characterized by 3 parameters ($a$, $b$, $k$). 
The SU$_q$(1,1) potential is also characterized by 3 parameters ($A$, $\tau$, 
$N$), while the WKB-EPs of the q-oscillator and Q-oscillator are characterized
 by 2 parameters ($\omega$, $\tau$ and $\omega$, $Q$ respectively) and the 
modified P\"oschl-Teller potential is characterized by only one parameter
($A$). It is therefore not surprising that a connection is possible only for
the case in which the number of parameters of the two potentials to be related
is the same.  

In summary, we have proved that the quasi-exactly soluble potentials 
corresponding to a harmonic oscillator with quartic and 
sextic anharmonicities  can be approximated by the WKB equivalent potentials 
corresponding to a deformed anharmonic oscillator with SU$_q$(1,1) 
symmetry. As a result one can use the relevant SU$_q$(1,1) oscillator in 
order to approximate the levels of the QESP which cannot be obtained using 
that approach. In the specific case under discussion
the extra levels obtained in this way are the levels of the same parity lying
above the first $n+1$ known levels, as well as all the levels of the 
opposite parity. The fact that the modified P\"oschl-Teller potential 
with SU(1,1) symmetry cannot be connected to these QESPs, while its 
deformed version with SU$_q$(1,1) symmetry can, is an example where the 
use of q-deformation allows the solution of an otherwise intractable problem. 

In ref. \cite{BAR} the deformed anharmonic oscillator with SU$_q$(1,1) 
symmetry 
has been proved appropriate for the description of vibrational spectra of 
diatomic molecules. An extension to vibrational spectra of highly symmetric 
polyatomic molecules has been given in \cite{BD3611}. The construction of 
QESPs 
related to realistic vibrational spectra is receiving attention.   

One of the authors (DB) has been supported by the EU under contract 
ERBCHBGCT930467. Another author
(HM) acknowledges the warm hospitality during his stay 
in NCSR ``Demokritos''.


\begin{thebibliography}{9}

\bibitem{PLB307}
D. Bonatsos and C. Daskaloyannis, Phys. Lett. B 307 (1993) 100. 

\bibitem{BAR}
D. Bonatsos, E. N. Argyres and P. P. Raychev, J. Phys. A 24 (1991) L403.

\bibitem{JMP33} D. Bonatsos, C. Daskaloyannis and K. Kokkotas, J. Math. Phys. 
33 (1992) 2958. 

\bibitem{TurUsh} 
A. V. Turbiner and A. G. Ushveridze, Phys. Lett. A 126 (1987) 181. 

\bibitem{Tur118} A. V. Turbiner, Commun. Math. Phys. 118 (1988) 467. 

\bibitem{Tur276} A. Turbiner, Phys. Lett. B 276 (1992) 95. 

\bibitem{Ioann}
 Th. Ioannidou, Diploma Thesis, U. Thessaloniki (1993), unpublished. 

\bibitem{SunLi} C. P. Sun and W. Li, Commun. Theor. Phys. 19 (1993) 191. 

\bibitem{BD3611}
 D. Bonatsos and C. Daskaloyannis, Phys. Rev. A 48 (1993) 3611. 

\end{thebibliography}
\end{document}